\documentclass[12pt,french]{amsart}
\usepackage{amssymb,amsmath,a4wide,babel}

\newtheorem{theorem}{Th\'eor\`eme}[section]
\newtheorem{lemma}[theorem]{Lemme}

\newcommand{\ep}{\varepsilon}

\newcommand{\pa}{\partial}
\newcommand{\RR}{{\mathbb{R}}}
\renewcommand{\div}{\mathop{{\rm div}}}

\title[P\'erim\`etre sur les vari\'et\'es]{P\'erim\`etre sur
les vari\'et\'es et application aux \'equations
aux d\'eriv\'ees partielles}%
\author{Satyanad Kichenassamy}%
\address{D. M. I., Ecole Normale Sup\'erieure, 45 rue d'Ulm,
75230 Paris Cedex 05,France}%

\thanks{S\'eminaire Equations aux D\'eriv\'ees Partielles 1986--1987,
Ecole Polytechnique, Expos\'e no.~14, 24 f\'evrier 1987.}%
\begin{document}
\maketitle

Le perim\`etre d'une partie mesurable de $\RR^N$ a \'et\'e
d\'efini dans \cite{DG} par E. de Giorgi ; c'est la variation
totale de sa fonction caract\'eristique. Il v\'erifie une
in\'egalit\'e isop\'erim\'etrique et permet d'\'etablir des
estimations \emph{a priori} dans de nombreux probl\`emes
lin\'eaires et non lin\'eaires (voir Talenti \cite{T} et ses
r\'ef\'erences).

Nous allons donner ici une g\'en\'eralisation de cette notion au
cas de parties de vari\'et\'es compactes.

Nous d\'efinirons dans une premi\`ere partie le p\'erim\`etre
d'une partie $E$ d'une vari\'et\'e riemannienne compacte
orient\'ee comme la variation totale de sa fonction
caract\'eristique $\chi_E$, et nous montrerons que ce
p\'erim\`etre est la limite des variations totales des
r\'egularis\'ees de $\chi_E$ par le noyau de la chaleur. Nous en
d\'eduirons une in\'egalit\'e isop\'erim\'etrique, et une formule
de type Fleming-Rishel.

Dans un second temps, nous \'etudierons un probl\`eme
quasilin\'eaire elliptique dans $\RR^N$, dont l'\'etude a
necessit\'e l'introduction du p\'erim\`etre sur les vari\'et\'es
(Kichenassamy \cite{K}). On verra que les m\'ethodes usuelles de
sym\'etrisation dans $\RR^N$ achoppent, mais que l'on pourra
conclure en introduisant une m\'ethode de sym\'etrisation sur
$S^N$; cette technique met en jeu les propri\'et\'es du
p\'erim\`etre evoqu\'ees dans la premi\`ere partie.

\section{P\'erim\`etre dans une vari\'et\'e.}
On consid\`ere dans cette section une vari\'et\'e riemannienne
$(M,g)$ compacte, orient\'ee et sans bord.

\subsection*{a) Premi\`ere d\'efinition}

Soit $u\in L^1(M)$. On d\'efinit sa \emph{variation totale}
\begin{equation}
 V(u) = \sup_{\scriptsize\begin{array}{l}
  \varphi \text{ champ de vecteurs }C^\infty\\
         \qquad\qquad|\varphi|\leq 1\end{array}} \int_M u \div \varphi\, dV .
\end{equation}
On d\'efinit le p\'erimetre d'une partie mesurable $E$ de $M$ par
\begin{equation}
P(E)=V(\chi_E)
\end{equation}
o\`u $\chi_E$ d\'esigne la fonction caract\'eristique de $E$.

\emph{Remarques} : i) Lorsque $E$ admet un bord de classe
$C^\infty$, $P(E)$ correspond \`a la mesure de $\pa E$ au sens
usuel (soit le volume de $\pa E$ pour la structure riemannienne
induite)

ii) $V(u) = V(1-u)$ pour tout $u$ et donc $P(E)=P(M\setminus E)$.

iii) Il est facile de voir que si $(u_n)$ tend vers $u$ dans
$L^1(M)$, alors $V(u)\leq\varliminf_{n\to\infty} V(u_n)$.

\subsection*{b) Deuxi\`eme d\'efinition}

Soit $u_0\in L^1(M)$. On pose $u(t)=\int_M
e_0(x,y,t)u_0(y)\,dV(y)$ o\`u $e_k(x,y,t)$ d\'esigne le noyau de
la chaleur pour les $k$-formes.
\begin{lemma} Il existe $\theta$ ne d\'ependant que de $M$ tel
que $f_{u_0} : t \mapsto e^{-\theta t}\int_M |du(t)|\,dV$ soit
d\'ecroissante.
\end{lemma}
On pose alors que le p\'erim\`etre de $E$ est la limite, lorsque
$t\downarrow 0$, de la fonction $f_{\chi_E}$. L'\'equivalence des
deux d\'efinitions r\'esulte du
\begin{lemma}
Pour tout $u\in L^1 (M)$,
\[
V(u) = \lim_{t\downarrow 0}\int_M |du(t)|\,dV
\]
o\`u $u$ est la solution de l'\'equation de la chaleur avec
donn\'ee initiale $u_0$.
\end{lemma}

\subsection*{Preuves de 1.1. et 1.2.} Nous donnons ici de br\`eves
indications sur ces d\'emonstrations, renvoyant \`a Kichenassamy
\cite{K} pour les d\'etails.

--- La preuve du Lemme 1.1. repose sur un calcul explicite
aboutissant \`a une relation de la forme
\[\frac d{dt} \int_M J(|du|)\,dV\leq C\int_M J(|du|)\,dV
\] o\`u $J$ est une fonction
convexe approchant la valeur absolue et $C$ est une borne sur les
termes de courbure intervenant dans les formules de Weitzenb\"ock.

--- Le lemme 1.2. repose essentiellement sur la propri\'et\'e
$d_x e_k = \delta_y e_{k+1}$ du noyau de la chaleur.

\subsection*{c) In\'egalite isop\'erim\'etrique}

\begin{theorem} Il existe une constante $C_I$ telle que pour tout
$E$ mesurable inclus dans $M$ on ait
\[ P(E)\geq C_I\min (|E|,|M-E|)^{1-1/N}.
\]
\end{theorem}
\emph{Preuve.} On d\'eduit de l'in\'egalit\'e de Sobolev pour les
r\'egularis\'ees de $u_0\in L^1(M)$ que
\[V(u_0) \geq \text{C}^\text{te}
\left(\inf_{c\in\RR} \int_M |u_0-c|^{\frac N{N-1}} \right)^{1-1/N}
\]
On applique ensuite ce r\'esultat \`a $u_0=\chi_E$.

\subsection*{d) Formule de type Fleming-Rishel}

Il s'agit du r\'esultat suivant :
\begin{theorem} Pour tout $u\in L^1(M)$,
$V(u)$ est fini si et seulement si $t\mapsto P(u> t)$ est dans
$L^1(\RR)$ et l'on a alors l'\'egalit\'e :
\[V(u) = \int_{-\infty}^\infty P(u>t)\,dt.
\]
\end{theorem}
\emph{Remarque} : Le r\'esultat correspondant dans le cas de
l'espace euclidien est d\^u \`a Fleming et Rishel \cite{FR}.

\subsection*{D\'emonstration} On montre d'abord que, pour tout
$u\in L^1(M)$, on a
\begin{equation}
  V(u) \leq \raisebox{-4pt}{${}_*$}\int_{-\infty}^\infty P(u>t)\,dt
\end{equation}
et
\begin{equation}
  V(u) \geq \raisebox{9pt}{${}^*$}\!\!\int_{-\infty}^\infty P(u>t)\,dt
\end{equation}
Le th\'eor\`eme en r\'esulte.

i) \emph{Preuve de (3) :} On \'ecrit, pour tout $x$ dans $M$,
\begin{equation}
  u(x) = \int_{-\infty}^\infty b(t,x)\,dt
\end{equation}
o\`u
\[  b(t,x)=
\begin{cases}
+1 & \text{ si }u(x)>t\geq 0 \\ -1 & \text{ si }u(x)\leq t <0\\ 0
& \text{dans tous les autres cas.}
\end{cases}
\]
Il est clair que $\int_M\int_{-\infty}^\infty |b(t,x)|\,dt\,
dV(x)=\|u\|_{L^1(M)}$.

Soit $\varphi$ un champ de vecteurs $C^\infty$. On a
\begin{equation} \int_M u \div \varphi
\;dV = \int_M dV(x)\,\int_{-\infty}^\infty b(t,x) \div \varphi\;dt
\end{equation}
et l'on voit facilement que $V(b(t,x)) = P(u>t)$ pour tout $t$.

Prenant $\varphi$ de longueur $\leq 1$, on obtient
\[
 \int_M u \div \varphi\;dV
 \leq\raisebox{-4pt}{${}_*$}\int_{-\infty}^\infty P(u>t)\,dt
\]
d'o\`u l'in\'egalit\'e (3).

ii) \emph{Preuve de (4) :}
Il r\'esulte de la deuxi\`eme d\'efinition du p\'erim\`etre
(Eq.~(3)) qu'il existe des fonctions $(u_n)_{n\geq 1}$ telles que
\begin{equation}
u_n\to u \text{ p.p.\ et dans }L^1(M) ;\quad
V (u_n ) \to V (u) ;\quad
u_n\in C^\infty (M).
\end{equation}
Par approximation, on peut supposer que les $u_n$ sont des
fonctions de Morse, auquel cas (par la formule de la "co-aire" par
exemple)
\begin{equation} V(u_n) = \int_{-\infty}^\infty P (u_n
> s) \,ds.
\end{equation}
Par le th\'eor\`eme de convergence domin\'ee, $\chi_{u_n >s}$ tend
vers $\chi_{u >s}$ dans $L^1(M)$, pour chaque $s$ fix\'e. Par la
Rem.~iii) du a) ci-dessus,
\[ \varliminf_{n\to\infty} P (u_n> s) \geq P (u> s).
\]
Le lemme de Fatou appliqu\'e \`a (8) donne l'in\'egalit\'e (4).

\emph{Remarque} : Il r\'esulte en particulier de ce th\'eor\`eme
que si $u\in W^{1,1}(M)$, $-\frac d{d t} \int_{u> t}|\nabla u| dV
= P (u> t)$ p.p.

\section{Application \`a un probl\`eme quasilin\'eaire}

\subsection*{a) Le probl\`eme}

On cherche $u$ solution du probl\`eme suivant :
\begin{equation}
  \left\{
  \begin{array}{c}
Au := -\div(|\nabla u|^{p-2}\nabla
u)=\sum_{i=1}^m\gamma_i\delta(x-a_i)\text{ dans }\RR^N\\ u\to
0\text{ lorsque } |x|\to\infty,
  \end{array}
  \right.
\end{equation}
o\`u $p>1$, $N\geq 2$, $\gamma_i\in\RR$, $a_i\in\RR^N$. On dira
que $u$ est \emph{solution} de (9) si $u$ est de classe $C^1$ dans $\RR^N
\setminus\{a_1,\dots,a_m\}$, $|\nabla u|^{p-1}\in L^1_{\text{loc}}(\RR^N)$
et si $u$ v\'erifie (9) au sens des distributions. Il est facile de
trouver une solution dans le cas $m=1$ (une seule singularit\'e) :
de fait, si $\varphi(x) = C_p|x|^{(p-N)/(p-1)}$ si $p\neq  N$
(resp.~$C_N$ Log$\,(1/|x|)$ pour $p = N$) avec $C_p$, $C_N$
convenablement choisis, alors $A\varphi = \delta$. Bien s\^ur,
$\varphi\to 0$ \`a l'infini ssi $p < N$. Nous allons maintenant
voir comment r\'esoudre (9) dans le cas g\'en\'eral :
\begin{theorem}
Supposons que $\sum_{i=1}^m\gamma_i = 0$ si $p \geq N$. Alors il
existe une seule solution $u$ de (9) telle que $u -
\sum_{i=1}^m\gamma_i\varphi(x-a_i)$ soit born\'ee sur $\RR^N$.
\end{theorem}

Pour la preuve de ce r\'esultat et les motivations de (9) on
renvoie le lecteur \`a \cite{K}. Nous donnerons ici bri\`evement
les \'etapes qui conduisent \`a l'existence d'une solution dans le
cas $p<N$ et dans le cas $p=N$, ce dernier utilisant le
p\'erim\`etre \'etudi\'e au \S\ 1.

\subsection*{b) Existence---$p < N$}

 On \'etablit ici
l'existence d'une fonction $u$ telle que $Au=0$ sur
$\RR^N\setminus
\{a_l,\dots,a_m\}$ et telle que $u -
\sum_{i=1}^m\gamma_i\varphi(.-a_i)$ soit born\'ee. On peut en
d\'eduire \cite{K} que l'\'equation (9) est v\'erifi\'ee et que
$u$ admet une limite lorsque $|x|\to\infty$, d'o\`u l'existence
d'une solution de notre probl\`eme.

\underline{1${}^{\text{e}}$ \'etape} : On r\'esout
\[\left\{
\begin{array}{rl}
Au^\ep = \chi^\ep & \text{ sur }B(1/\ep)\\ u^\ep = 0& \text{ sur
}\pa B(1/\ep)
\end{array} \right.
\]
pour tout $\ep>0$, en prenant pour $\chi^\ep$ des fonctions
approchant la mesure $\sum_{i}\gamma_i\delta(x-a_i)$
(convenablement choisies).

\underline{2${}^{\text{e}}$ \'etape} : On montre, par une
technique de sym\'etrisation, que pour tout compact $K$ de
$\RR^N$, il existe $C_K$ tel que
\[\|u^\ep\|_{L^q(K)}\leq C_K.
\]
Ici, $q$ est un r\'eel $> p-1$ (ind\'ependant de $K$).

\underline{3${}^{\text{e}}$ \'etape} : On montre, gr\^ace \`a des
th\'eor\`emes de r\'egularit\'e pour l'\'equation $Au = 0$ (voir \cite{K})
qu'il existe un $\alpha$ positif tel que pour tout compact $K$ de
$\RR^N\setminus\{a_1,\dots,a_m\}$,
\[\|u^\ep\|_{C^{1,\alpha}(K)}\leq C'_K.
\]
Il est alors possible de passer \`a la limite.

La 2${}^{\text{e}}$ \'etape utilise une majoration des
sym\'etris\'ees $(u^\ep)^*$ par un multiple de la fonction
$\varphi$ ; elle est possible parce que $\varphi\to 0$ pour
$|x|\to\infty$. Lorsque $p = N$, il faut un autre argument.

\subsection*{c) Existence---$p=N$}

Il suffit d'amender la 2${}^{\text{e}}$ \'etape du b). Pour ce
faire, on se ram\`ene, gr\^ace a l'invariance conforme du
probl\`eme, \`a borner les solutions de classe $W^{1,N}$ de
\[ -\div{}_{S^N}(|\nabla u|^{N-2}\nabla u)=f\text{ sur la sph\`ere
}S^N, \text{ o\`u }\|f\|_{L^1}\leq \text{C}^\text{te}.
\]
et ce en termes de la seule quantit\'e $\|f\|_{L^1}$. Posons
$\mu(t)= \text{mes}\,(u> t)$. Normalisons $u$ (par addition d'une
constante) de sorte que
\begin{equation}
\text{mes}\, (u > 0) \text{ et mes }\, (u < 0) \text{ soient }\leq
\frac12  \text{mes}\,(S^N).
\end{equation}
On va maintenant estimer s\'epar\'ement les parties positive et
n\'egative de $u$.

Soit \underline{$t> 0$}. (10) et l'in\'egalit\'e
isop\'erim\'etrique donnent
\begin{equation}
P(u>t) \geq C_I \mu(t)^{1-1/N}.
\end{equation}
D'autre part, pour presque tout $t$ on a
\begin{equation}
-\frac{d}{dt}\int_{u>t}|\nabla u| \,dV=P(u>t).
\end{equation}

Multipliant l'\'equation satisfaite par $u$ par $(u-t)^+$, on
obtient
\begin{equation}
-\frac{d}{dt}\int_{u>t}|\nabla u|^N \,dV\leq \|f\|_{L^1}
\end{equation}
pour presque tout $t$.

Enfin, $t\mapsto t^{-1/(N-1)}$ \'etant convexe, on a (p.p en t)
\begin{equation}
  \left(
  \frac{\frac{d}{dt}\int_{u>t}|\nabla u|^N dV}
       {\frac{d}{dt}\int_{u>t}|\nabla u| dV}
  \right)^{-\frac 1{N-1}}
  \leq \frac{\mu'}{\frac{d}{dt}\int_{u>t}|\nabla u|dV}.
\end{equation}
Combinant (11)-(14) on obtient $-\mu'\geq C\mu$ et comme $\mu$ est
non-croissante, on en d\'eduit que la sym\'etris\'ee
(uni-dimensionnelle) de $u^+$ est major\'ee par une fonction de la
forme $s \mapsto a-b\,\text{Log\,}s$. Il en r\'esulte une borne
$L^q$ pour $u^+$, pour \underline{tout $q> 1$}. On estime de
m\^eme la fonction $u^-$. On ach\`eve la d\'emonstration comme
dans le cas $p < N$.

\end{document}